\documentclass{amsproc}

\usepackage{amsthm}

\usepackage{graphicx}
\usepackage{epsfig}
\usepackage{epstopdf}

\usepackage{latexsym}
\usepackage{amsmath}
\usepackage{amssymb}
\usepackage{amsfonts}

\theoremstyle{definition}

\newtheorem{theorem}{Theorem}%[section]
\newtheorem{lemma}[theorem]{Lemma}

\theoremstyle{definition}
\newtheorem{definition}[theorem]{Definition}

\newtheorem{corollary}[theorem]{Corollary}
\newtheorem{prop}[theorem]{Proposition}

\theoremstyle{remark}
\newtheorem{remark}[theorem]{Remark}

\def\proof{\par\medskip\noindent{\it Proof. }}

\def\R{{\mathbb R}}
\def\C{{\mathbb C}}
\def\Z{{\mathbb Z}}

\def\N{{\mathbb N}}

\def\al{\alpha}
\def\be{\beta}
\def\ga{\gamma}

\def\la{\lambda}

\def\om{\omega}

\newcommand\s{{}^{*}{\R}}
\newcommand\n{{}^{*}{\N}}

\newcommand\us{{}^{u}{\R}}
\newcommand\un{{}^{u}{\N}}
\renewcommand\c{\mathfrak c}

\begin{document}

\title{Krull dimensions of rings of holomorphic functions}
\author{Michael Kapovich}
\address{Department of Mathematics, 
University of California, Davis, 
CA 95616}
\email{kapovich@math.ucdavis.edu}
\thanks{The author was supported in part by the NSF Grant DMS-12-05312 and by the Korea Institute for Advanced Study (KIAS)}



\date{December 27, 2015}

\begin{abstract}
We prove that the Krull dimension of the ring of holomorphic functions of a connected complex manifold is at least the cardinality of continuum iff it is $>0$. 
\end{abstract}

\maketitle

Let $R$ be a commutative ring. Recall that the {\em Krull dimension} $\dim(R)$ of $R$ is the supremum of cardinalities lengths of chains of  distinct proper prime ideals in $R$. Our main result is: 

\begin{theorem}\label{main}
Let $M$ be a connected complex manifold and $H(M)$ be the ring of holomorphic functions on $M$. Then the Krull dimension of $H(M)$ either equals $0$ (iff $H(M)= \C$) or is infinite, iff $M$ admits a nonconstant holomorphic function $M\to \C$. More precisely, unless $H(M)=\C$, $\dim H(M)\ge \c$, i.e., the ring $H(M)$ 
contains a chain of distinct prime ideals whose length has cardinality of continuum. 
\end{theorem}

Our proof of this theorem mostly follows the lines of the proof by Sasane \cite{Sasane}, who proved that for each nonempty domain $M\subset \C$ the Krull dimension of $H(M)$ is infinite (he did not prove that $\dim H(M)\ge \c$). 

\begin{remark}
We note that Henricksen \cite{Henricksen} was the first to prove that the Krull dimension of the ring of entire functions on $\C$ has cardinality at least continuum.
\end{remark}

In our proof we will use the Axiom of Choice  in two ways: (a) to establish existence of certain maximal ideals and  
(b) to get existence of a nonprincipal ultrafilter $\om$ on $\N$ and, hence of the ordered field $\s$ of 
{\em nonstandard real} (or, {\em surreal}) numbers. The field $\s$ contains $\n$, the {\em nonstandard natural} (or {\em surnatural}) numbers.

The field $\s$ is a certain quotient of the countable direct product $\prod_{k\in \N} \R$; we will denote the equivalence class (in $\s$) of a sequence $(x_k)$ in $\R$  by $[x_k]$. Accordingly, $\n$ consists of equivalence classes $[n_k]$ of sequences of natural numbers. Roughly speaking, we will use $\n$ and certain order relation on it to compare rates of growth of sequences of  natural numbers. 

\begin{definition}\label{defn:ample}
A commutative unital ring $R$ is {\em ample} if there exists a sequence of valuations $\nu_k$ on $R$ such that for each $\be\in \n$, there $a=a_\be\in R$ with the property
\begin{equation}\label{1}
[\nu_{k}(a)]= \be  . 
\end{equation}
\end{definition}

\medskip
The main technical result of this paper is: 

\begin{theorem}\label{thm:T}
For each ample ring  $R$, $\dim(R)\ge \c$.  In particular, $R$ has infinite Krull dimension. 
\end{theorem}

This theorem and its proof are inspired by Theorem 2.2 of \cite{Sasane}, although some parts of the proof resemble the ones of \cite{Henricksen}. 

We will verify, furthermore, that whenever $M$ is a connected complex manifold which has a nonconstant holomorphic function, the ring $H(M)$ is ample.  This, combined with Theorem \ref{thm:T}, will immediately imply Theorem \ref{main}.

\begin{remark}
1. We refer the reader to Section 5.3 of \cite{Clark} for further discussion of algebraic properties of rings of holomorphic functions.

2. Theorem \ref{main} shows that for every Stein manifold $M$ (of positive dimension), the ring $H(M)$ has infinite Krull dimension. In particular, this applies to any noncompact connected Riemann surfaces (since every such surface is Stein, \cite{BS}). 

3. Noncompact connected complex manifolds $M$ of dimension $>1$ can have  $H(M)=\C$; for instance, take $M$ to be the complement to a finite subset in a compact connected complex manifold (of dimension $>1$). 
\end{remark}

\medskip
\noindent {\bf Acknowledgements.} This note grew out of the mathoverflow question, 

 http://mathoverflow.net/questions/94537, 

\noindent and I am grateful to Georges Elencwajg for asking the question. I am also grateful to Pete Clark for pointing at several errors in earlier versions of the paper (most importantly, pointing out that Lemma \ref{lem:max} is needed for the proof, which forced me to use ultrafilters) and providing references.

\section{Surreal numbers}
 
We refer the reader to \cite{Goldblatt} for a detailed treatment of surreal numbers, below is a brief introduction. 
A nonprincipal ultrafilter on $\N$ can be regarded as a finitely-additive probability measure on $\N$ which vanishes on each finite subset and takes  the value $0$ or $1$ on each subset of $\N$. The existence of nonprincipal ultrafilters (the {\em ultrafilter lemma}) follows from the Axiom of Choice. Subsets of full measure are called {\em $\om$-large}. Using $\om$ one defines the following equivalence relation on the product 
$$
\prod_{k\in \R} \R. 
$$
Two sequences $(x_k)$ and $(y_k)$ are equivalent if $x_k=y_k$ for an $\om$-all $k$, i.e. the set 
$$
\{k: x_k=y_k\}
$$ 
is $\om$-large. The quotient by this equivalence relation, denoted  
$$
\s=  \prod_{k\in \N} \R/\om, 
$$
is the set of surreal numbers. Let $[x_k]$ be the equivalence class of the sequence $(x_k)$. 

The binary operations on sequences of real numbers project to binary operations on $\s$ making $\s$ a field. The total order  $\le$ on $\s$ is defined by $[x_k]\le [y_k]$ iff $x_k\le y_k$ for an $\om$-all $k\in \N$. With this order, $\s$ becomes an ordered field.  

The set of real numbers embeds into $\s$ as the set of equivalence classes of constant sequences; the image of a real number $x$ under this embedding is still denoted $x$. We set $\s_+:=\{\al\in \s: \al>0\}$. 

The projection of
$$
\prod_{k\in \N} \N  \subset \prod_{k\in \N} \R
$$
to $\s$ is denoted $\n$, this is the set of {\em surnatural numbers}. We define a further equivalence relation $\sim_u$ on 
$\s$ by:
$$
\al\sim_u \be
$$
if there exist positive real numbers $a, b$ such that
$$
a\al \le \be \le b \al. 
$$
The equivalence class $(\al)$ of $\al\in \s$ (for this equivalence relation) is a multiplicative analogue of the {\em galaxy} $gal(\al)$ of $\al$, see   \cite{Goldblatt}:

\begin{definition}
The {\em galaxy} $gal(\al)$ of a surreal number $\al\in \s$ is the union
$$
\bigcup_{n\in \N} [\al-n, \al+n] \subset \s. 
$$ 
In other words, $\be\in gal(\al)$ iff there exist a real number $a$ such that $\al-a\le \be \le \al+ a$. 
\end{definition}

The next lemma is immediate:

\begin{lemma}
For $\al\in \s_+$, the equivalence class $(\al)$ of $\al$ equals $\exp( gal ( \log(\al)))$. 
\end{lemma}

We let $\us$ denote the quotient $\s/\sim_u$ and $\un$ the projection of $\n$ to $\us$. Define the total order 
$\gg$ on $\us$ by
$$
(\be) \gg (\al)  
$$  
if for every real number $c$, $c\al < \be$. By abusing the notation, we will simply say that $\be \gg \al$, with $\al, \be\in \s$. 

For the reader who prefers to think in terms of sequences of (positive) real numbers, 
the relation $(\be) \gg (\al)$ is an analogue of the relation  
$$
(a_n)= o((b_n)), \quad n\to \infty. 
$$

\begin{remark}
The equivalence relation $\sim_u$ and the order $\gg$ are similar to the ones used by Henricksen in \cite{Henricksen}. \end{remark}

%The following proposition is most likely already in the literature. 

\begin{prop}\label{prop:un}
The set $\un$ has the cardinality of continuum. 
\end{prop}
\proof Note first, that $\s$ has cardinality of continuum, hence, the cardinality of $\un$ is at most $\c$. The proof of the proposition then reduces to two lemmata.

\begin{lemma}\label{lem:cont}
The set $gal(\s_+)$ of galaxies $\{gal(\al): \al\in \s_+\}$ has the cardinality of continuum.
\end{lemma}
\proof For each $\al=[a_k]\in \s_+$, the  galaxy $gal(\al)$ contains the surnatural number $\lceil \al \rceil= [ b_k]$, 
where $b_k=   \lceil a_k \rceil$. For each surnatural number $\be\in \n$, and natural number $n\in \N$, the intersection
$$
[\be -n, \be+ n]\cap \n
$$
is finite, equal $\{\be -n,..., \be + n\}$. Therefore, $gal(\be) \cap \n= \{\be\} + \Z$. It follows that the map 
$$
\n \to gal(\s_+), \quad \be \mapsto gal(\be) 
$$
is a bijection modulo $\Z$. Lastly, the set of surnatural numbers $\n$ has the cardinality of continuum. \qed 

\begin{lemma}\label{lem:surj}
The map $\lambda: \n\to gal(\s_+)$, $\lambda: \be\mapsto gal(\log(n))$,  is surjective. 
\end{lemma}
\proof For each $\al\in \s_+$ let $\be = \lceil \exp(\al) \rceil \in \n$. Since $\log(x+1) -\log(x) \le 1$ for $x\ge 1$, we have that
$$
\log (\be) \in  gal(\al). \qed 
$$ 

Now, we can finish the proof of the proposition. The map $\la: \n\to gal(\s_+)$ descends to a map $\mu: \un\to gal(\s_+)$. According to Lemma \ref{lem:surj}, the map $\mu$ is surjective. By Lemma \ref{lem:cont}  the set  $gal(\s_+)$ has the cardinality of continuum.  \qed 

\medskip 
We will prove Theorem \ref{thm:T} in the next section by showing that for each ample ring $R$, the ordered set 
$(\un, \gg)$ embeds into the poset of prime ideals in $R$ reversing the order:
$$
(\be) \gg (\al) \Rightarrow P_\beta\subsetneq P_\al
$$
for certain prime ideals $P_\ga\subset R$ determines by $(\ga)\in \un$.  
Proposition \ref{prop:un} will then imply that the Krull dimension of $R$ is at least $\c$.

\section{Krull dimension of ample rings}

\noindent Recall  that a {\em valuation} on a unital ring $R$ is a map $\nu: R\to \R_+\cup \{\infty\}$ such that:

1. $\nu(a+b)\ge \min(a, b)$,

2. $\nu(ab)=\nu(a)+\nu(b)$. 

3. $\nu(a)=\infty \iff a=0$. 

4. $\nu(1)=0$.

\noindent For the following lemma, see  Theorem 10.2.6 in \cite{Cohn} (see also Proposition 4.8 of \cite{Clark} or Theorem 1 in \cite{Kaplansky}).

\begin{lemma}\label{L0}
Let $I$ be an ideal in a commutative ring $A$ and $M\subset A\setminus I$ be a subset closed under multiplication. Then there exists an ideal $J\subset A$ containing  $I$ and disjoint from $M$, so that $J$ is maximal with respect to this property. Furthermore, $J$ is a prime ideal in $A$. 
\end{lemma}

Let $R$ be an ample ring and $\nu_k$ the corresponding sequence of valuations on $R$. 
For each $\be\in \n$ we define %the ideal 
$$
I_\be:=\{a\in R | ~~ [\nu_k(a)] \gg [\be] \}\subset R. 
$$  

\begin{lemma}\label{lem:add}
Each $I_\al$ is an ideal in $R$. 
\end{lemma}
\proof We will check that $I_\al$ is additive since it is clearly closed under multiplication by elements of $R$. 
Take $p', p''\in I_\al$,  
$$
[\nu_k(p')]\gg \al, [\nu_k(p'')] \gg \al. 
$$
By the definition of a valuation,  
$$
n_k:=\nu_k(p'+ p'') \ge \min( \nu_k(p'), \nu_k(p'')),
$$
for each $k\in \N$. For $m\in \N$, define the $\om$-large sets 
$$
A'=\{k: \nu_k(p')\ge m \al\}, \quad A''=\{k: \nu_k(p'')\ge m \al\}.$$ 
Therefore, their intersection $A=A'\cap A''$ is $\om$-large as well, which implies that 
$$
\forall m\in \N, [n_k] \ge m \al \Rightarrow [n_k] \gg \al.  \qed 
$$

%($I_\be$ is clearly closed under multiplication by elements of $R$, for additivity of $I_\be$ 
Then for each $\ga \gg \be$, the element $a_\ga$ as in Definition \ref{defn:ample}, belongs to $I_\beta$.
It follows that $I_\be\ne 0$ for every $\be$. Define the subsets    
$$
M_\be:= \{ a\in R  | \exists n\in \N, [\nu_k(a)] \le n \be \}\subset R; 
$$
each $M_\be$ is  closed under the multiplication. It is immediate that whenever $\al\le \be$, we have the inclusions 
$$
I_\be \subset I_\al , \quad M_\al\subset M_\be. 
$$
It is also clear that $I_\be \cap M_\be=\emptyset$. At the same time,  for each 
$\be \gg \al$,
$$
a_\beta\in I_\al \cap M_\beta. 
$$

 For each $\al$ we let 
${\mathcal J}_\al$ denote the set of ideals $P\subset R$ such that 
$$
I_\al\subset P, P\cap  M_\al=\emptyset.$$ 
By Lemma \ref{L0}, every maximal element $P\in {\mathcal J}_\al$ is a prime ideal. 

\begin{lemma}\label{lem:max}
Every ${\mathcal J}_\al$ contains unique maximal element, which we will denote $P_\al$ in what follows. 
\end{lemma}
\proof Suppose that $P', P''$ are two  maximal elements of  ${\mathcal J}_\al$. We define the ideal 
$P=P' + P''$. Clearly, $P$ contains $I_\al$. To prove that $P$ is disjoint from $M_\al$, take  
 $p'\in P', p''\in P''$, since $p'\notin M_\al, p''\notin M_\al$. Then the same proof as in Lemma \ref{lem:add} shows that 
 $[\nu_k(p'+ p'') ]\gg \al$ which means that $p'+p''\notin M_\al$. 
 Thus, $P\in {\mathcal J}_\al$ and, in view of maximality of $P', P''$, we obtain
$$
P'= P= P''. \qed 
$$

For each $\be \gg \al$ we define the ideal $Q_{\al\be}:=I_\al+ P_{\be}$.

\begin{lemma}\label{L1}
$Q_{\al\be} \cap M_\al=\emptyset$. 
\end{lemma}
\proof The proof is similar to the one of the previous lemma. 
Let $q=c+p$, $c\in I_\al, p\in P_{\be}$. Since $p\notin M_{\be}$, $p\notin M_\al$ as well. Therefore, 
$$
[\nu_{k}(p)] \gg \al.   
$$ 
Since $c\in I_\al$, 
$$
[\nu_{k}(c)] \gg \al. 
$$
Hence, 
$$
[\nu_{k}(c+p)] \gg \al 
$$
as well. Thus, $q\notin M_\al$. \qed 

\begin{corollary}\label{cor:C}
$Q_{\al\be}\in {\mathcal J}_\al$. In particular, $Q_\al\subset P_\al$. 
\end{corollary}
\proof It suffices to note that $I_\al\subset Q_{\al\be}$ according to the definition of $Q_{\al\be}$. \qed 

\begin{lemma}\label{L2}
The inequality $\be \gg \al$ implies $P_{\be}\subset P_\al$ and this inclusion is proper. 
\end{lemma}
\proof By the definition of $Q_{\al\be}$ and Corollary \ref{cor:C}, we have the inclusions
$$
P_{\be}\subset Q_\al\subset P_\al. 
$$
We now claim that $P_{\be}\ne Q_{\al\be}=I_\al+ P_{\be}$. Recall that $a_\al\in I_\al\subset Q_{\al\be}$ and 
$a_\al\in M_{\be}$, while $M_{\be}\cap P_{\be}=\emptyset$. Thus, $a_\al\in Q_{\al\be} \setminus P_{\be}$. \qed

\medskip
According to Proposition \ref{prop:un}, the set  $\n$ of surnatural numbers contains a subset $S$ of cardinality continuum such that for all $\al< \be$ in $S$, we have
$\be \gg \al$. The map
$$
\al\mapsto P_\al
$$
sends each $\al\in S$ to a prime ideal in $R$; $\al< \be$ implies that $P_\be \subsetneq P_\al$. 

We conclude that the ring $R$ contains the (descending) chain of distinct prime ideals $P_\al, \al\in S$; the length of this chain has the cardinality of continuum.
In particular, $\dim(R)\ge \c$. Theorem \ref{thm:T} follows. \qed

\section{Ampleness of rings of holomorphic functions}

\noindent We will need the following classical result, see e.g. \cite[Ch. VII, Theorem 5.15]{Conway}: 

\begin{theorem}\label{thm:W}
Let $D\subset \C$ be a domain, and let $c_k\in D$ be a sequence which does not accumulate anywhere in $D$ and let $m_k$ be a sequence of natural numbers. Then there exists a holomorphic function $g$ in $D$ which has zeroes only at the points $c_k$ and such that $m_k$ is the order of zero of $g$ at $c_k$, $k\in \N$. 
\end{theorem}

\begin{corollary}\label{cor:main}
If $M$ is a connected complex manifold which admits a nonconstant holomorphic function $h: M\to \C$, 
then the ring $H(M)$ is ample. 
\end{corollary}
\proof We let $D$ denote the image of $h$. Pick a sequence $c_k\in D$ which converges to a point in 
$\hat{\C} \setminus  D$ and which consists of regular values of $h$.  
(Here $\hat{\C}$ is the Riemann sphere.) For each $c_k$ the preimage $C_k:= h^{-1}(c_k)$ 
is a complex submanifold in $M$; in each $C_k$ pick a point $b_k$. Define valuations 
$$\nu_k: H(M)\to \Z_+ \cup \{\infty\}$$ 
by $\nu_k(f):= ord_{b_k}(f)$, the total order of $f$ at $b_k$, cf.  \cite[Chapter C, Definition 1]{Gunning}.  

Now, given $\be\in \n$, $\be=[m_k]$, we let $g=g_\be$ denote a holomorphic function on $D$ as in Theorem 
\ref{thm:W}.   Define $a=a_\be:= g\circ h\in H(M)$. Then 
$\nu_k(a)= m_k$, which implies that the ring $H(M)$ is ample. \qed 

\medskip 
Ampleness of $H(M)$ together with Theorem \ref{thm:T} imply Theorem \ref{main}.

%\newpage


\begin{thebibliography}{BLP05}

\bibitem[BS]{BS}
H. Behnke, K. Stein, {\em 
Entwicklung analytischer Funktionen auf Riemannschen Fl\"{a}chen}, 
Math. Ann. Vol. {\bf 120} (1949) p. 430--461. 

\bibitem[Cla]{Clark}
P. Clark, ``Commutative Algebra''. Preprint. %, http://math.uga.edu/~pete/integral.pdf. 

\bibitem[Coh]{Cohn}
P. M. Cohn, ``Basic Algebra: Groups, Rings and Fields'', Springer Verlag,  2004.  

\bibitem[Con]{Conway}
J. B. Conway, ``Functions of One Complex Variable'' I, 2nd edition. Springer Verlag, New York (1978). 

%\bibitem[DK]{DK}
%C. Drutu, M. Kapovich, ``Geometric Group Theory'',\newline http://www.math.ucdavis.edu/EPR/ggt.pdf, to appear in  the AMS Colloquium series. 

\bibitem[Go]{Goldblatt}
 R. Goldblatt, ``Lectures on the hyperreals,''
Graduate Texts in Mathematics, Vol.  {\bf 188}, Springer-Verlag, 1998. 

%\bibitem[Go]{G}
%M. Golsi\'nski, {\em Henriksen's contributions to residue class rings of analytic and entire
%functions}, Topology and its Applications, Vol. {\bf 158} (2011) p. 1756--1761. 

\bibitem[Gu]{Gunning}
R. Gunning, ``Introduction to Holomorphic Functions of Several Variables,'' Volume {\bf 1}, Wadsworth \& Brooks/Cole, 1990. 

\bibitem[H]{Henricksen} 
M. Henricksen, {\em On the prime ideals of the ring of entire functions}, Pacific J. Math. Vol. {\bf 3} (1953) p. 711--720.


\bibitem[K]{Kaplansky}
I. Kaplansky, ``Commutative Rings'', Allyn and Bacon, Inc., Boston, Mass. 1970. 

\bibitem[S]{Sasane}
A.\ Sasane, 
{\em On the Krull dimension of rings of transfer functions}, 
Acta Appl. Math. Vol. {\bf 103} (2008) p. 161--168. 



\end{thebibliography}
\end{document}